\documentclass[a4paper]{ifacconf}

\usepackage{graphicx,amsmath,url}      
\usepackage[round]{natbib}             

\def\portugues{1} 

\ifx\portugues\undefined
\def\portugues{0}
\fi

\if\portugues0
   \usepackage[english]{babel}
  \else
   \usepackage[spanish,brazil,english]{babel}
\fi

\usepackage[T1]{fontenc}

\usepackage[utf8]{inputenc}

\usepackage{ae}

\if\portugues1
 { }
 { }
 { }
 
\fi
            
\begin{document}

\if\portugues1


%
\selectlanguage{brazil}
 
\begin{frontmatter}

\title{Cálculo do maior expoente de Lyapunov usando estimador recursivo com fator variável} 


\author[]{Emily A. de Sousa} 
\author[]{Erivelton G. Nepomuceno} 	
\author[]{Márcio F. S. Barroso}

\address[]{GCOM - Grupo de Controle e Modelagem \\ Departamento de Engenharia Elétrica\\ Universidade Federal de São João del-Rei\\ Pça. Frei Orlando, 170 - Centro - 36307-352 - São João del-Rei, MG, Brasil, (e-mail: emilyasousa10@gmail.com, nepomueno@ufsj.edu.br, barroso@ufsj.edu.br)}

\selectlanguage{english}
\renewcommand{\abstractname}{{\bf Abstract:~}}
\begin{abstract}                
Chaotic systems have been investigated in the most diverse areas. One of the first steps in chaotic system research is the detection of chaos. The largest Lyapunov exponent (LLE) is one of the most widely used techniques for this purpose. Recently, techniques for calculating LLE have been developed taking into account the error due to the finite precision of computers. Recursive methods were employed to improve such algorithms. However, this method uniformly weighed the data used to calculate the LLE. This paper investigates the different weighing of the data based on the hypothesis that the initial data have a higher precision than the last data processed by the algorithm, since the process is subject to error accumulation. In five tested systems, the proposed method obtained more accurate results in three. However, there was an increase in the variance of the result obtained for two of the evaluated systems.

\vskip 1mm
\selectlanguage{brazil}
{\noindent \bf Resumo}: Sistemas caóticos têm sido investigados nas mais diversas áreas. Uma das primeiras etapas na pesquisa em sistema caóticos é a detecção da presença de caos. O maior expoente de Lyapunov (LLE) é uma das técnicas mais empregadas para esse fim. Recentemente, técnicas para o cálculo do LLE tem sido desenvolvidas levando em conta o erro devido a precisão finita dos computadores. Métodos recursivos foram empregados para aprimorar tais algoritmos. Entretanto, tal  método ponderou de modo uniforme os dados utilizados para o cálculo do LLE. Este trabalho investiga a ponderação diferenciada dos dados pautada na hipótese de que os dados iniciais possuem uma maior precisão do que os últimos dados processados pelo algoritmo, uma vez que o processo está sujeito ao acúmulo do erro. Em cinco sistemas testados, o método proposto obteve resultados mais precisos em três. Entretanto, observou-se um aumento da variância do resultado obtido para dois dos sistemas avaliados. 
\end{abstract}

\selectlanguage{english}


\selectlanguage{brazil}

\end{frontmatter}
\else
%

\begin{frontmatter}

\title{Style for SBA Conferences \& Symposia: Use Title Case for
  Paper Title\thanksref{footnoteinfo}} 

\thanks[footnoteinfo]{Sponsor and financial support acknowledgment
goes here. Paper titles should be written in uppercase and lowercase
letters, not all uppercase.}

\author[First]{First A. Author} 
\author[Second]{Second B. Author, Jr.} 
\author[Third]{Third C. Author}

\address[First]{Faculdade de Engenharia Elétrica, Universidade do Triângulo, MG, (e-mail: autor1@faceg@univt.br).}
\address[Second]{Faculdade de Engenharia de Controle \& Automação, Universidade do Futuro, RJ (e-mail: autor2@feca.unifutu.rj)}
\address[Third]{Electrical Engineering Department, 
   Seoul National University, Seoul, Korea, (e-mail: author3@snu.ac.kr)}
   
\renewcommand{\abstractname}{{\bf Abstract:~}}   
   
\begin{abstract}                
These instructions give you guidelines for preparing papers for IFAC
technical meetings. Please use this document as a template to prepare
your manuscript. For submission guidelines, follow instructions on
paper submission system as well as the event website.
\end{abstract}


\end{frontmatter}
\fi


\section{Introdução}

Sistemas caóticos são investigados nas mais diversas áreas. Novas técnicas de controle foram desenvolvidas para tratar das especificidades desta categoria de sistemas não-lineares \citep{OGY1990} e grande interesse tem sido dedicado para o sincronismo de osciladores caóticos \citep{PCJ+1997}. Além disso, inúmeras aplicações têm sido observadas na literatura \citep{torres2005,souza2012,cuomo1993,Yeh2012}, sendo que mais recentemente um grande número de trabalhos tem sido publicados na área de criptografia \citep{NNAB2019}.

Uma das primeiras abordagens na investigação de sistemas caóticos é o discernimento se o sistema é efetivamente caótico ou não. O cálculo do maior expoente de Lyapunov (LLE) é uma das técnicas mais empregadas no problema de detecção da presença de caos \citep{Rosenstein1993,PEIXOTO201836}. Além disso, o LLE pode ser explorado para detectar e controlar órbitas periódicas instáveis \citep{Perc2004} e adquirir informações importantes sobre a dinâmica do sistema. Os expoentes de Lyapunov quantificam a divergência exponencial das trajetórias em espaço de estado inicialmente próximas e estimam a caoticidade em um sistema. Nesse sentido, vários métodos numéricos para estimar o expoente de Lyapunov são encontrados na literatura \citep{Wolf1985,Bochi2008a,Matsuoka2015a,Kantz1994}.

Recentemente, \cite{Nepomuceno2016} implementaram um algoritmo simples, robusto, rápido e fácil para estimar o LLE. O método proposto é baseado no conceito do limite inferior do erro, introduzido por \cite{NM2016}. Esse método proposto caracteriza-se por não precisar de nenhum tipo de parametrização, dimensão de incorporação e estimativa de fluxo linearizado, e usa a equação dinâmica original. No entanto, há limitações como  a dificuldade de elaborar uma extensão intervalar natural e a complexa manipulação em algumas equações. Nesse sentido, a fim de suprir os problemas encontrados em \cite{Nepomuceno2016}, \cite{PEIXOTO201836} desenvolveram um trabalho no qual o modo de arredondamento, de acordo com o padrão IEEE 754-2008, foi aplicado em vez de utilizar extensões intervalares naturais para calcular o LLE. Nesse trabalho foi possível observar que os dois modos de arredondamento diferentes apresentam efeitos semelhantes quando produzidos por duas extensões intervalares. Dessa forma, usando o estimador recursivo de mínimos quadrados, foi alcançado vantagens no método proposto, obtendo sucesso na aplicação em cinco exemplos clássicos de sistemas caóticos.

Entretanto, nos resultados atingidos por \cite{PEIXOTO201836}, as diferenças de precisão proveniente da implementação finita do cálculo não foram levadas em conta \cite{Nep2014,NMAR2017}. Uma vez que os dados utilizados para o cálculo são computados com precisão finita, pode-se levantar dúvidas se a recursividade utilizada em \cite{PEIXOTO201836} deveria assumir ponderações distintas ao longo de seu cálculo. Nesse sentido, esse trabalho objetiva investigar a perda de informações a cada interação. Para tanto, será proposto neste trabalho o uso de um fator $\lambda$, que pondera os dados de maneira diferenciada. Mais especificamente, é proposto que a ponderação varia inversamente com o número de iterações.  Até então, os estimadores recursivos utilizados em \cite{PEIXOTO201836} pesam igualmente os dados recentes e os dados do passado. Neste trabalho, sugere-se que o fator $\lambda$ tenha maior ponderação para os primeiros dados processados pelo algoritmo. Acredita-se que devido a acumulação de erros dos computadores, os dados subsequentes tenham menor precisão e assim influenciam negativamente no cálculo do LLE que realizam os procedimentos para o cálculo do LLE. Em cinco sistemas testados, o método proposto obteve resultados mais precisos em três. Entretanto, observou-se um aumento da variância para dois dos sistemas avaliados.

O trabalho foi dividido da seguinte forma: A Seção 2 apresenta os Conceitos Preliminares, que mostra um resumo das principais definições para o melhor entendimento da proposta do trabalho. Na Seção 3 é apresentado o Método Proposto. A seção 4 mostra os resultados obtidos, bem como comparações dos valores numéricos encontrados na literatura. Por fim, as conclusões estão na Seção 5.

\section{Conceitos Preliminares}

\subsection{Modo de arredondamento de acordo com o padrão IEEE 754-2008}
Desde as décadas de 1970 e 1980 o Institute of Electrical and Eletronics Engineers \citep{Ove2001} normalizou alguns aspectos para fabricantes de computadores e consultores de \textit{software} como a representação consistente em todas as máquinas que adotam o padrão e operação de arredondamento adequado para o ponto flutuante. Em 2008 uma versão atualizada \citep{IEE2008} foi publicada, com o objetivo de facilitar o movimento de programas existentes de diversos computadores, aprimorar as capacidades e a segurança disponíveis para usuários e programadores, fornecer suportes e proporcionar desenvolvimento de algumas funções elementares, aritmética de alta precisão, computação algébrica numérica e simbólica e habilitar o refinamento de extensões adicionais.

Alguns resultados de operações com ponto flutuante não podem ser representados exatamente, em razão da representação de números fora do intervalo normalizado ou quando o número de bits não é suficiente para representar o número com exatidão, dessa forma, é preciso aproximar o número desejado \citep{Gol1991,IEE2008,Ove2001}. O padrão IEEE define os modos de arredondamentos, que mapeiam um número real estendido para um número de ponto flutuante incluído nesse formato. Se \textit{x} é um número ponto flutuante então \textit{round(x) = x}, senão, o valor depende do modo de arredondamento. Definido \textit{round(x) = x\textsubscript{inf}} como o arredondamento para baixo, \textit{round(x) = x\textsubscript{sup}} como o arredondamento para cima, \textit{round(x) = x\textsubscript{inf}(x \textgreater 0)} ou \textit{round(x) = x\textsubscript{sup}(x \textless 0)} arredondamento em direção a zero, e em arredondamento para o mais próximo \textit{round(x)} é tanto x\textsubscript{inf} ou x\textsubscript{sup}, dependendo de qual for mais próximo de \textit{x}. Se houver um empate, aquele com o último bit significativo igual a zero é escolhido .

Utilizando o \textit{software Matlab}, é possível escolher se o número será arredondado para o menor, maior ou o ponto flutuante mais próximo \citep{Rum2005}, por meio de um comando interno: \verb|system_dependent(`setround',mod)|, onde \verb|mod=-Inf| altera o modo de arredondamento para baixo, \verb|mod=+Inf| para cima e  \verb|mod=0.5| o modo de arredondamento é definido como o mais próximo.

\subsection{Cálculo do LLE}
O expoente de Lyapunov está relacionado à convergência ou divergência exponencialmente rápida de órbitas próximas no espaço de fase. Como as órbitas próximas correspondem a estados quase idênticos, a divergência orbital exponencial significa que o sistema comportará de maneira bem diferente do esperado, assim, a capacidade de previsão é rapidamente perdida, caracterizando, dessa forma, o caos \citep{Wolf1985}. 

Mapas unidirecionais, como estudados nesse trabalho, possuem um único expoente de Lyapunov. Sistemas com esse expoente igual a zero significa uma órbita marginalmente estável, expoente negativo retrata órbitas periódicas, enquanto expoente positivo caracteriza um sistema caótico. Dessa forma, o cálculo do expoente de Lyapunov é explorado, e métodos para esse fim são propostos por vários autores \citep{Wolf1985,Kantz1994,Rosenstein1993}.

Recentemente alguns métodos foram propostos, a fim de simplificar, facilitar e dar maior robustez. \cite{Nepomuceno2016} propuseram um método baseado no limite inferior do erro (LBE), introduzido pela primeira vez em \cite{NM2016}, no qual é necessário somente as equações dinâmicas do sistema sob investigação. \cite{PEIXOTO201836} também propôs um método, utilizando o modo de arredondamento de acordo com o padrão IEEE 754-2008, que, assim como o trabalho de \cite{Nepomuceno2016}, utiliza mínimos quadrados para ajustar uma linha de inclinação do valor absoluto do erro do algoritmo natural definido, a inclinação da linha é o maior expoente positivo de Lyapunov.

\begin{table*}[!ht]
	\centering
		\caption{\label{tab:1} Sistemas caóticos explorados. A segunda e a terceira colunas mostram as descrições dos sistemas, enquanto a ultima coluna mostra o intervalo das condições iniciais adotadas em cada sistema para obter a média dos LLEs .}
	\begin{tabular}{lccc}
		\hline\hline \\[0.01cm]
		Sistema $\qquad$ & Equações                 & Parâmetros & Intervalo das condições iniciais \\[0.2cm] \hline
		Logístico         & $x_{n+1}=\mu x_n(1-x_n)$ & $\mu=4,0$  &  0,1 a 1      \\[0.2cm]
		H\'enon         & $x_{n+1}=1-ax_n^2+y_n$     & $a=1,4$   & -0,9 a 1       \\[0.2cm]
		&$y_{k+1}=bx_{k}$    & $b=0,3 $    &\\[0.2cm]
		Mapa Seno       & $x_{n+1}=ax_n - bx_n^3$  & $a=2,6868, b=0,2462$   & -3,8 a 3,8      \\[0.2cm]
		Mapa Tenda      & $x_{n+1}=r\min\{x_n,1-x_n\}$   & r=1,99  & 0,1 a 0,99      \\[0.2cm]
		Mackey-Glass    & $\dot{x} = \cfrac{ax_\tau}{1-x_\tau^c} -bx$ & $a=0,2$,  $b=0,1$,   $c=10 $, $  \tau = 30 $  &  0,15 a 1,5 	\\[0.2cm]
		\hline\hline	
		
	\end{tabular}
\end{table*}

O método proposto neste trabalho tem como base o algoritmo proposto por \cite{PEIXOTO201836}:
\begin{enumerate}
\item Escolha dois modos de arredondamento, tais como o arredondamento para o mais próximo e arredondamento para cima;
\item Com o mesmo software, hardware, sistema operacional, condições iniciais  e esquema de discretização, simule o sistema com dois modos de arredondamento previamente escolhidos;
\item Use o algoritmo recursivo de mínimos quadrados para estimar a inclinação do valor absoluto do algoritmo natural do LBE. A inclinação dessa linha é o LLE.
\end{enumerate}

\subsection{Estimador Recursivo de Mínimos Quadrados com parâmetros variantes no tempo}
As técnicas recursivas para a estimação de parâmetros são muito úteis \citep{Silva2008}. Quando os dados são medidos e disponibilizados sequencialmente, ou seja, a cada período de amostragem, um sistema qualquer de coleta de dados fornece medições correspondentes aquele instante é necessário o uso de estimadores recursivos. Além disso, esses estimadores são vantajosos na solução de problemas numéricos cuja solução em batelada seria difícil \citep{aguirre}.

Uma maneira muito utilizada para determinar um parâmetro desconhecido de uma equação de condições é minimizando a soma dos quadrados dos resíduos. O método dos mínimos quadrados, introduzido pela primeira vez por \cite{gauss} é utilizado na determinação do estimador recursivo de mínimos quadrados (RMQ), apresentado pela Equação \eqref{eq1} \citep{aguirre}.

\begin{equation}
\left\lbrace
\begin{array}{l}
K_k=\frac{P_{k-1}\mathbf{\psi}_k}{\mathbf{\psi}^T_k P_{k-1}\mathbf{\psi}_k+1},\\
\hat{\mathbf{\theta}}_k=\hat{\mathbf{\theta}}_{k-1}+K_k\left[y(k)-\mathbf{\psi}^T_k\hat{\mathbf{\theta}}_{k-1}\right],\\
P_k=P_{k-1}-K_k\mathbf{\psi}^T_kP_{k-1},
\end{array}
\right.
\label{eq1}
\end{equation}
em que $K$ é a matriz de ganho, $P$ é a matriz de covariância, $\psi$ é a matriz dos regressores e  $y(k)$ é o vetor da variável dependente.

A Equação \eqref{eq1} pondera de forma idêntica as observações obtidas ao longo do tempo. Contudo, os sistemas analisados neste trabalho variam no tempo, assim, algumas observações devem ser mais influentes na estimação dos parâmetros. Dessa forma, o RMQ pode ser redefinido, agora com um fator $\lambda$, que será responsável em ponderar de maneira diferenciada os dados observados. Então, o estimador recursivo de mínimos quadrados com o fator $\lambda$ é mostrado na Equação \eqref{eq2}.

\begin{equation}
\left\lbrace
\begin{array}{l}
K_k=\frac{P_{k-1}\mathbf{\psi}_k}{\mathbf{\psi}^T_k P_{k-1}\mathbf{\psi}_k+ \lambda},\\
\hat{\mathbf{\theta}}_k=\hat{\mathbf{\theta}}_{k-1}+K_k\left[y(k)-\mathbf{\psi}^T_k\hat{\mathbf{\theta}}_{k-1}\right],\\
P_k=\frac{1}{\lambda(k)}(P_{k-1}-\frac{P_{k-1}\mathbf{\psi}_k\mathbf{\psi}^T_kP_{k-1}}{\mathbf{\psi}^T_k P_{k-1}\mathbf{\psi}_k+ \lambda(k)}).
\end{array}
\right.
\label{eq2}
\end{equation}
\section{Metodologia}
O método proposto neste trabalho para o cálculo do LLE foi baseado nos métodos descritos em \cite{PEIXOTO201836} e \cite{Nepomuceno2016}. Neste caso, será usado o estimador recursivo de mínimos quadrados com parâmetros variantes no tempo, empregando o fator $\lambda$.

O fator $\lambda$ varia conforme os dados são processados. Os dados recebidos primeiro terão maior peso, isso em razão do armazenamento e propagação de erros do computador e da possível perda de informações do sistema, portanto, os dados no princípio são mais confiáveis. Neste trabalho $\lambda$ varia de acordo com cada novo conjunto de dados. O método de variação do fator $\lambda$ não compromete o comportamento do algoritmo \citep{barros,beza}. Neste trabalho, $\lambda(k)$ foi ajustado heuristicamente como $\lambda(k) = 1 + 1,02/k$, significando uma diminuição deste valor à medida que o $k$ cresce.

O LLE foi calculado em cinco importantes sistemas dinâmicos caóticos. Como os valores são bem sensíveis as condições iniciais o LLE tem uma certa variabilidade. Dessa forma, a fim de melhorar os resultados, o LLE foi calculado pela média de cem condições iniciais, definidas observando por onde os dados dos sistemas variam.

\section{Resultados}
\begin{table*}[!ht]
    \centering
    \caption{\label{tab:2} Cálculo do LLE em diferentes métodos para comparação dos resultados. Os valores esperados são os indicados na segunda coluna, que são os obtidos na literatura. As colunas 3 e 4 são referentes a média e desvio padrão quando o fator $\lambda$ não é empregado no estimador recursivo de mínimos quadrados. Já as colunas 5 e 6 indicam a média e o desvio padrão com o fato $\lambda$. Esses resultados evidenciam o aprimoramento do cálculo do maior expoente de Lyapunov com a técnica proposta.}
    \begin{tabular}{lccccc}
    	\hline\hline
    	                &  \\[0.01cm]
    	Sistema $\quad$ & Literatura &        $\mu$        &      $\sigma$       &        $\mu$        &      $\sigma$       \\
    	                &            & sem fator $\lambda$ & sem fator $\lambda$ & com fator $\lambda$ & com fator $\lambda$ \\[0.1cm] \hline
    	Logístico       &   0,6930   &       0,7014       &       0,1212       &       0,7098       &       0,0756       \\[0.1cm]
    	H\'enon         &   0,4180   &       0,4088       &       0,0947       &       0,4174       &       0,0847       \\[0.1cm]
    	Mapa Seno       &   0,7730   &       0,7487       &       0,1069       &       0,7598       &       0,0777       \\[0.1cm]
    	Mapa Tenda      &   0,6880   &       0,6898       &       0,0051       &       0,6891       &       0,0054       \\[0.1cm]
    	Mackey-Glass    &   0,0074   &       0,0074       &       0,0020       &       0,0082        &       0,0023        \\ \hline\hline
    \end{tabular}
\end{table*}

Os resultados apresentados foram obtidos usando o Matlab R2017a em um computador intel core i3 e com precisão dupla (64 bits). O método proposto foi aplicado nos sistemas caóticos descritos na Tabela \ref{tab:1}.

O comportamento a longo prazo de uma condição inicial especificada com qualquer incerteza não pode ser previsto, dessa forma, os valores numéricos para o LLE são muitos sensíveis às condições iniciais. Para melhorar os valores obtidos foi calculado o LLE para cem condições iniciais, então, foi obtido a média e o desvio padrão, para, assim, verificar o quanto esses valores estão próximos dos valores encontrados na literatura.

A Tabela \ref{tab:2} mostra a comparação entre os valores da literatura e os resultados obtidos utilizando o modo de arredondamento e o estimador recursivo de mínimos quadrados. Para a comparação foi feito a média e o desvio padrão, com 100 condições iniciais, quando os dados do sistema são igualmente ponderados, indicados nas colunas 3 e 4, e quando a perda de informações e a imprecisão dos cálculos foram considerados, indicados nas colunas 5 e 6. 

\begin{figure*}[htbp]
	\begin{center}
		\begin{tabular}{cc}           
			(a)                         &       (f)          \\
			\includegraphics[width=0.28\linewidth]{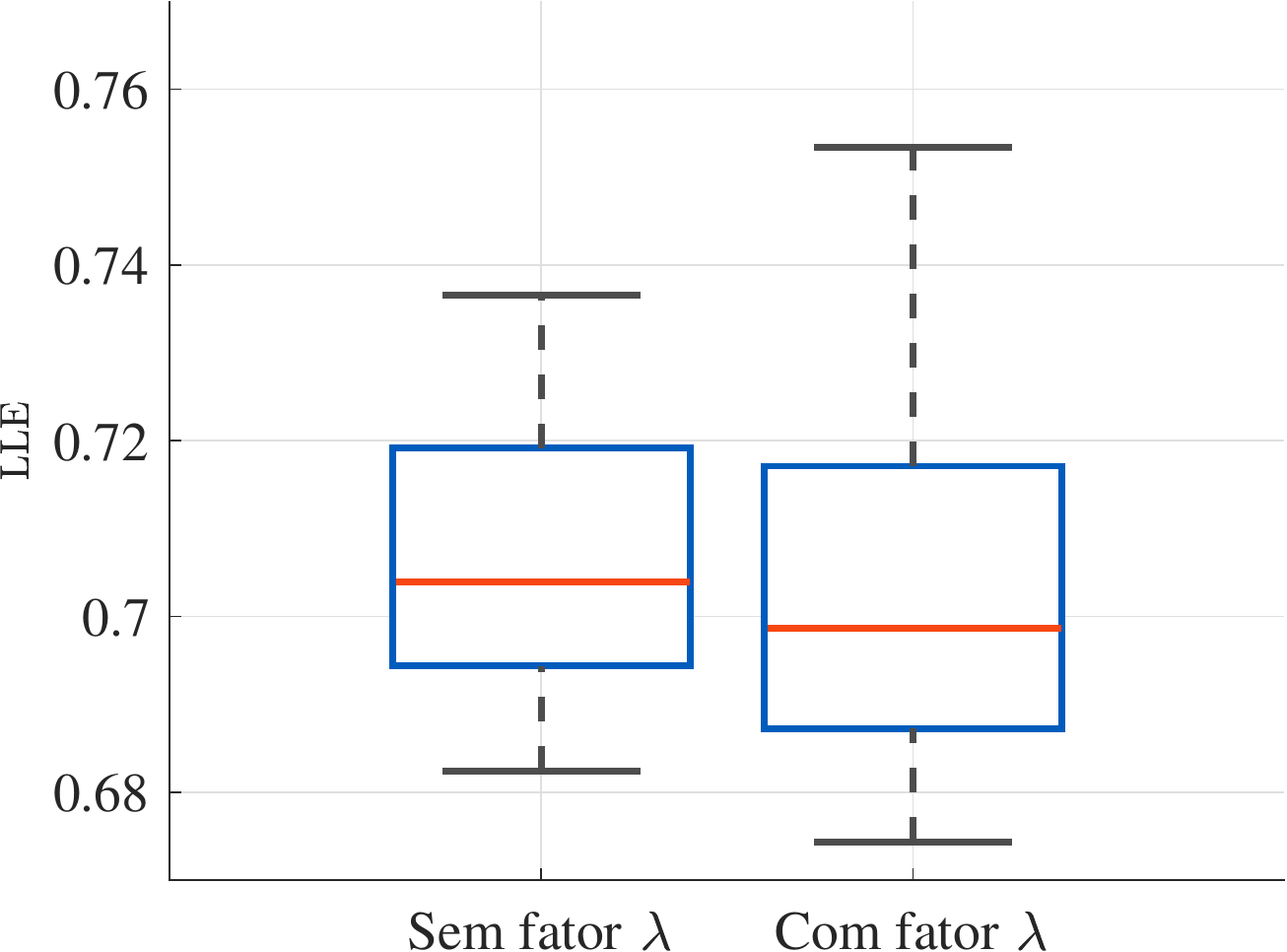} &
			\includegraphics[width=0.28\linewidth]{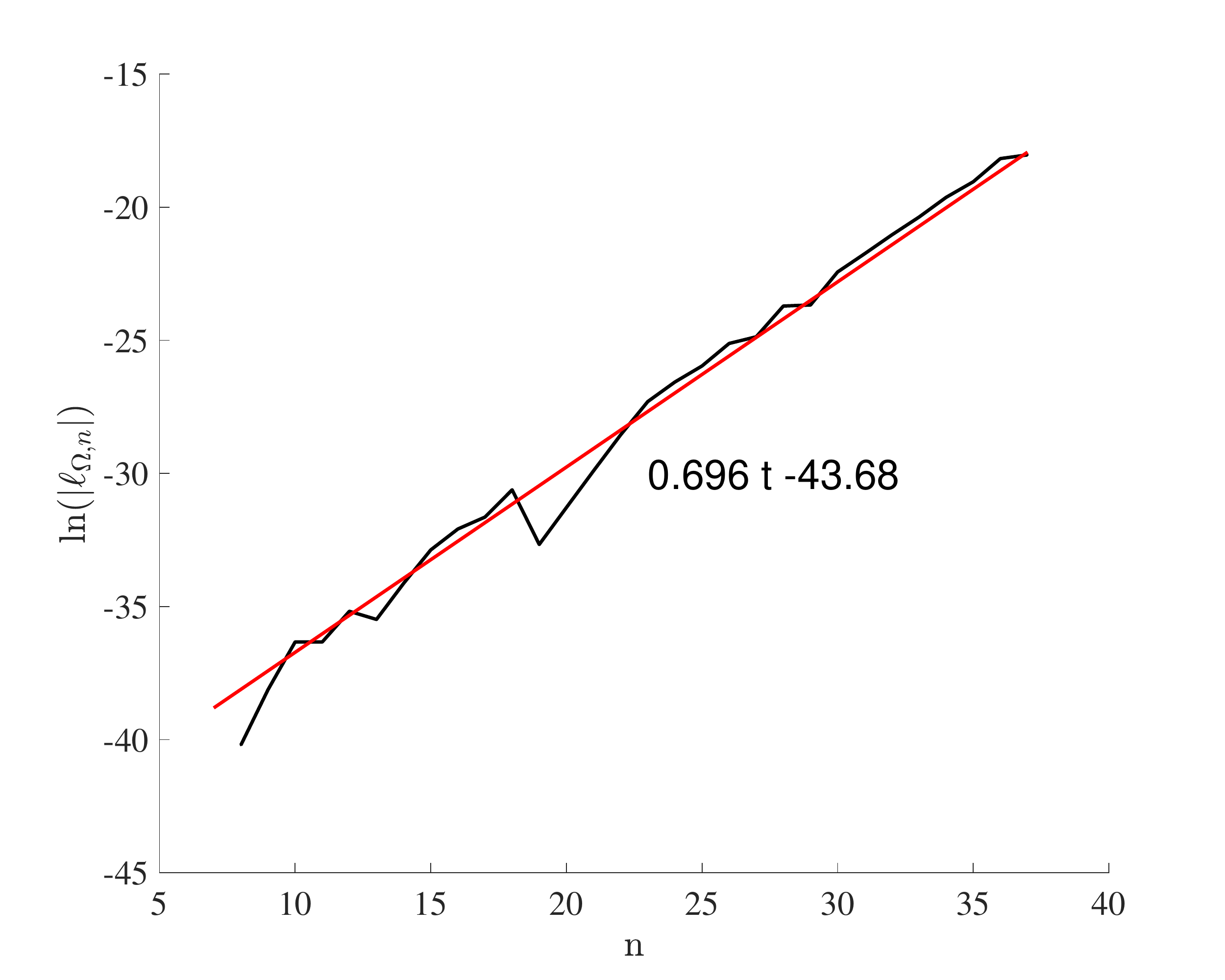}\\
 			(b) & (g) \\
 			\includegraphics[width=0.28\linewidth]{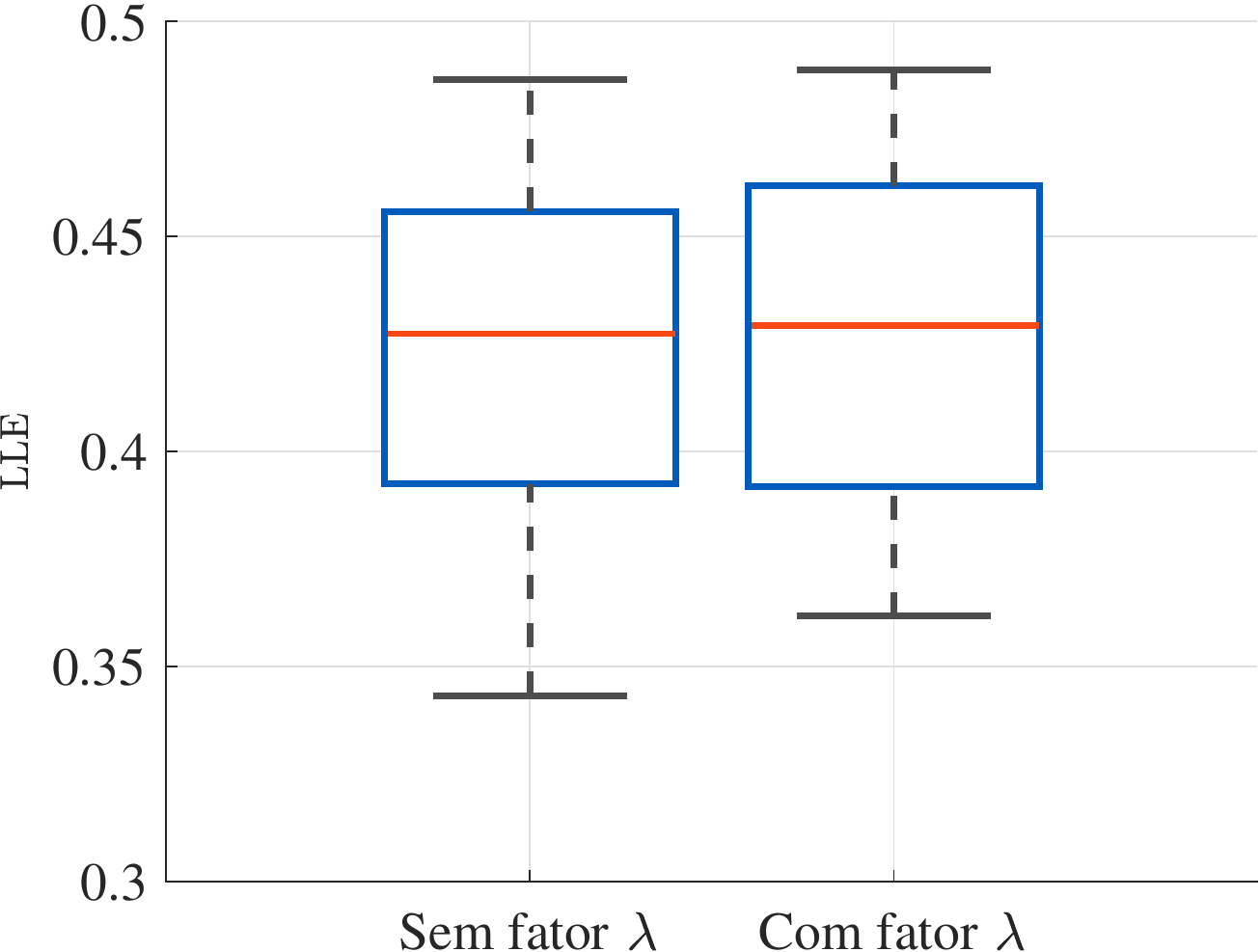}&
 			\includegraphics[width=0.28\linewidth]{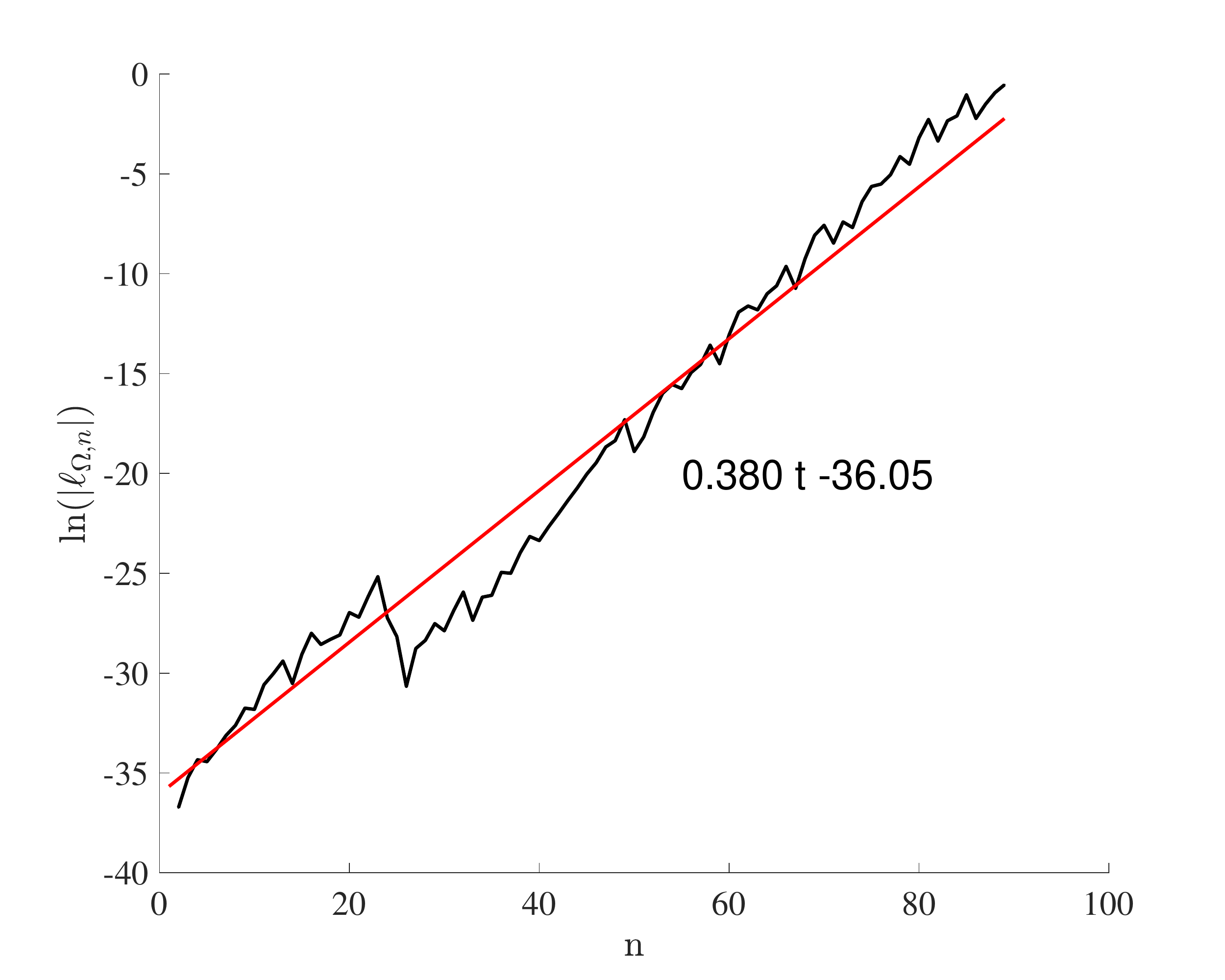} \\
 			(c)                         &           (h)       \\
     			\includegraphics[width=0.28\linewidth]{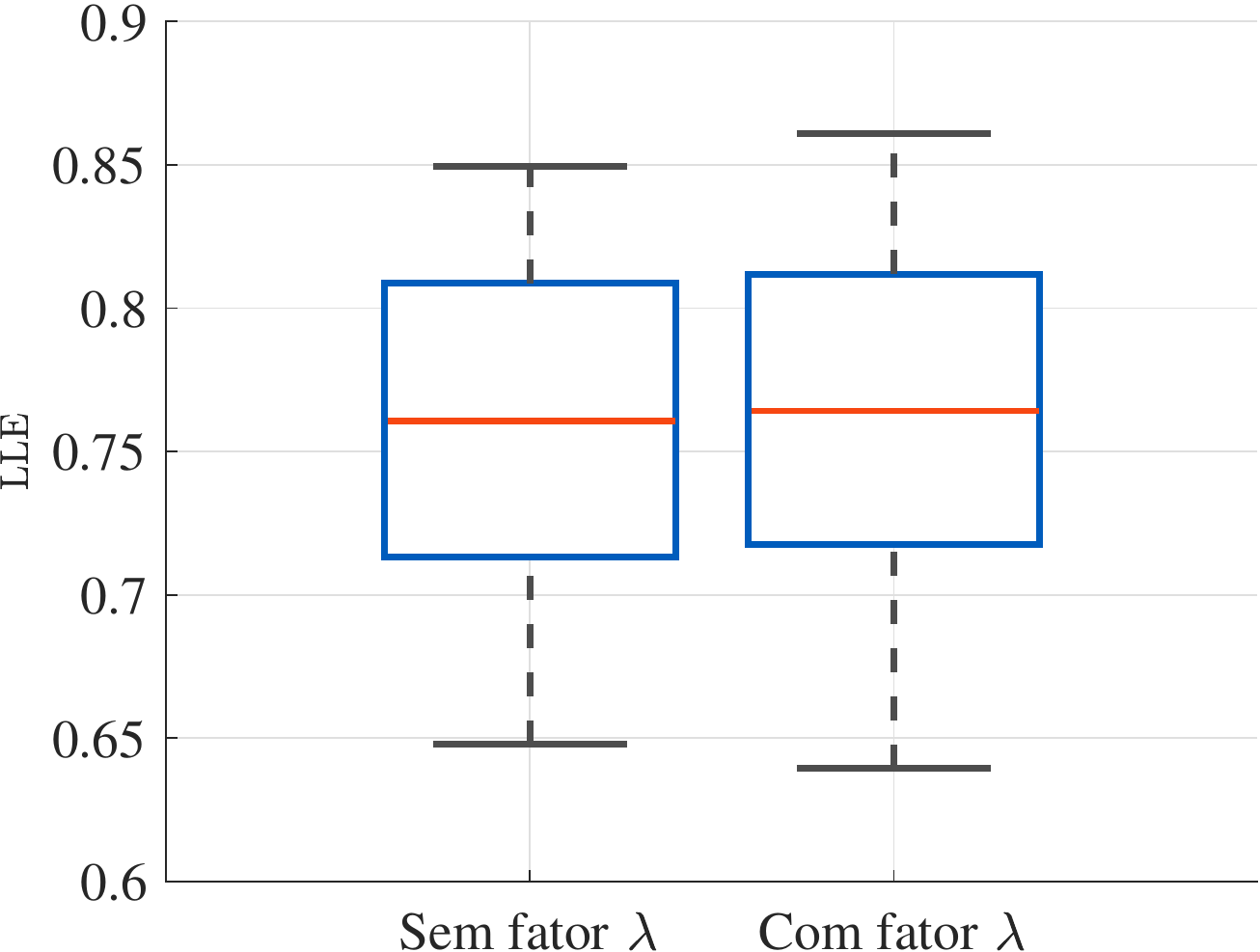}& 
     			\includegraphics[width=0.28\linewidth]{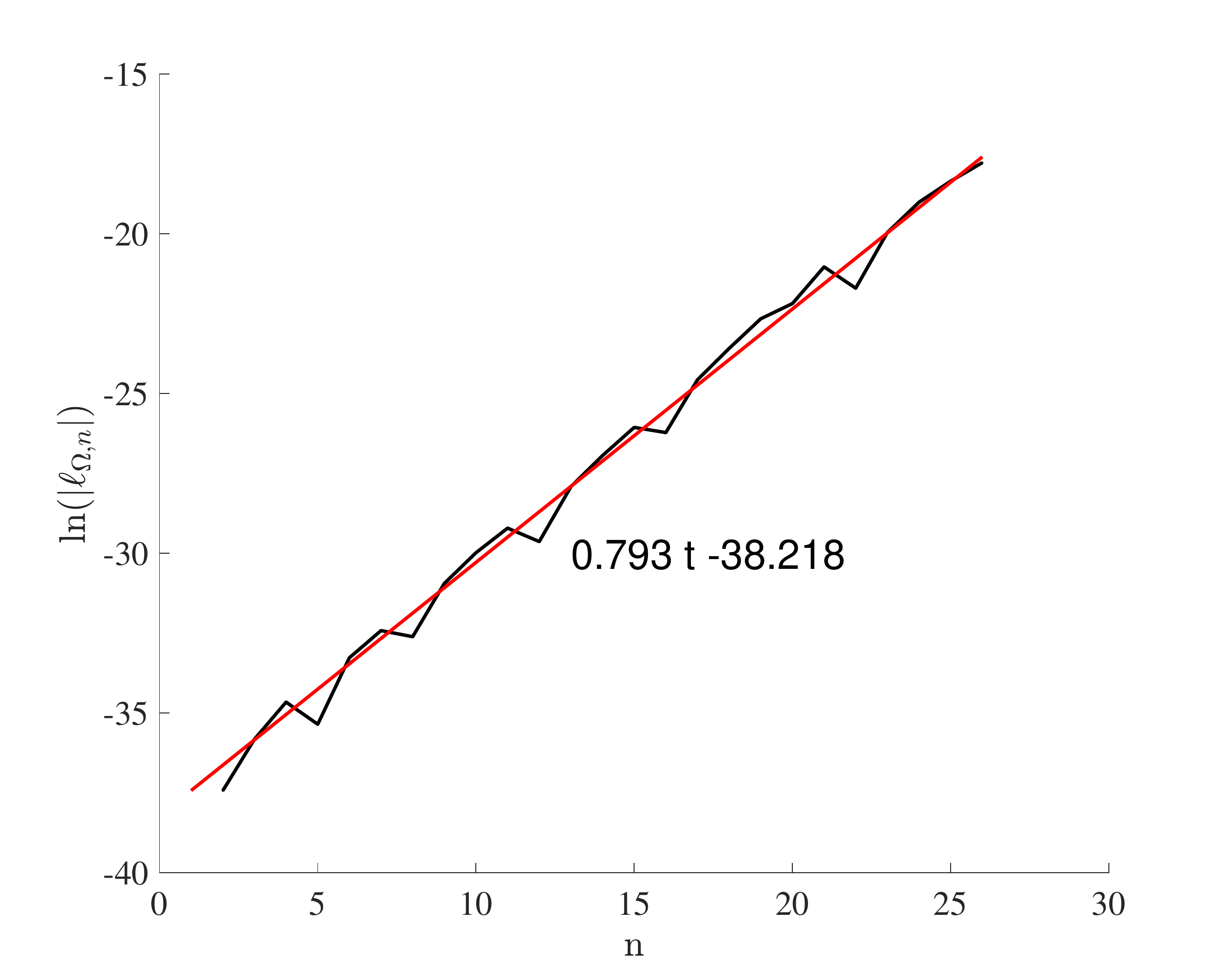}	\\
			
 			(d) & (i)\\
 			\includegraphics[width=0.28\linewidth]{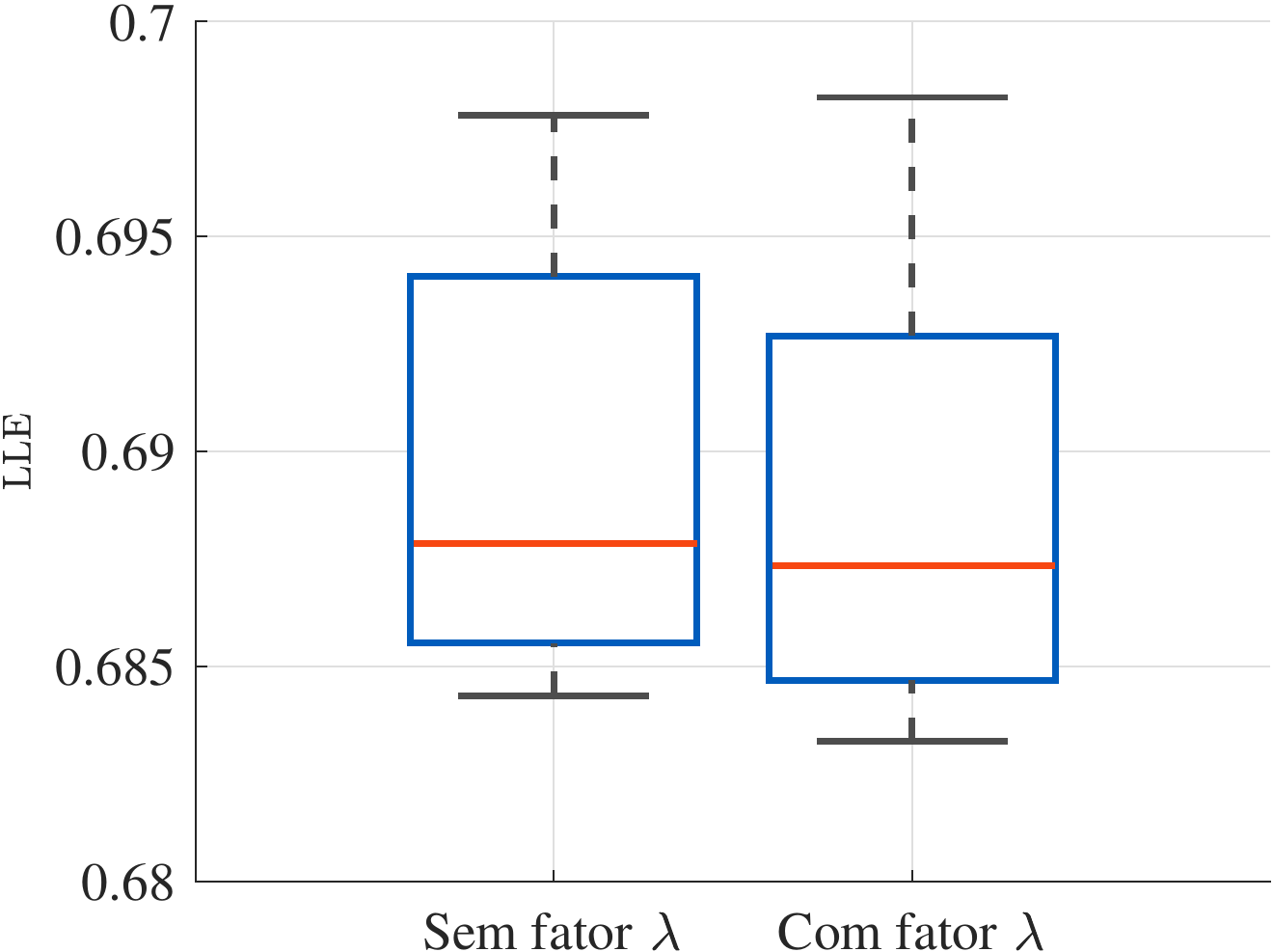} &
 			\includegraphics[width=0.28\linewidth]{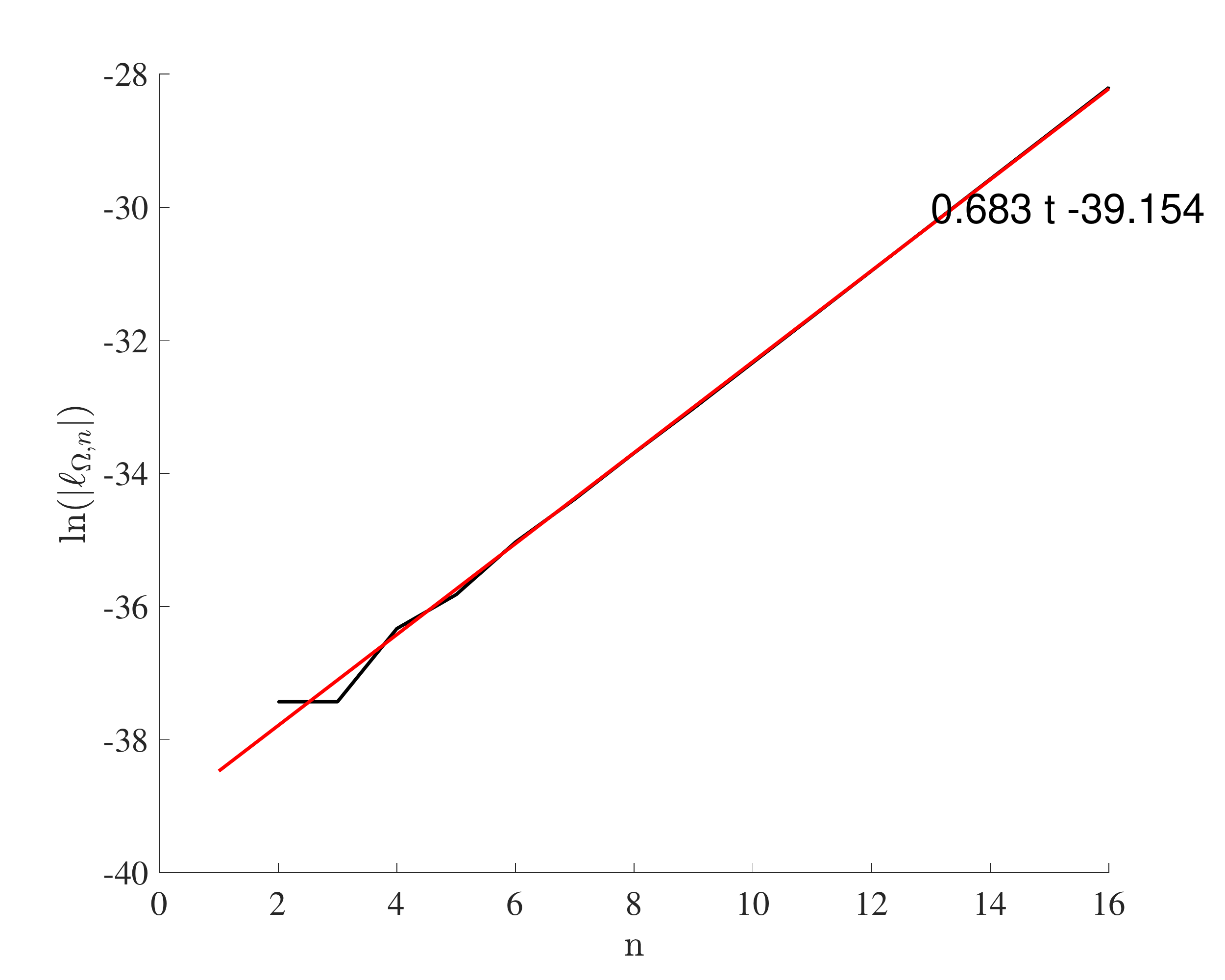} \\
 			(e)            & (j)                        \\
 			\includegraphics[width=0.28\linewidth]{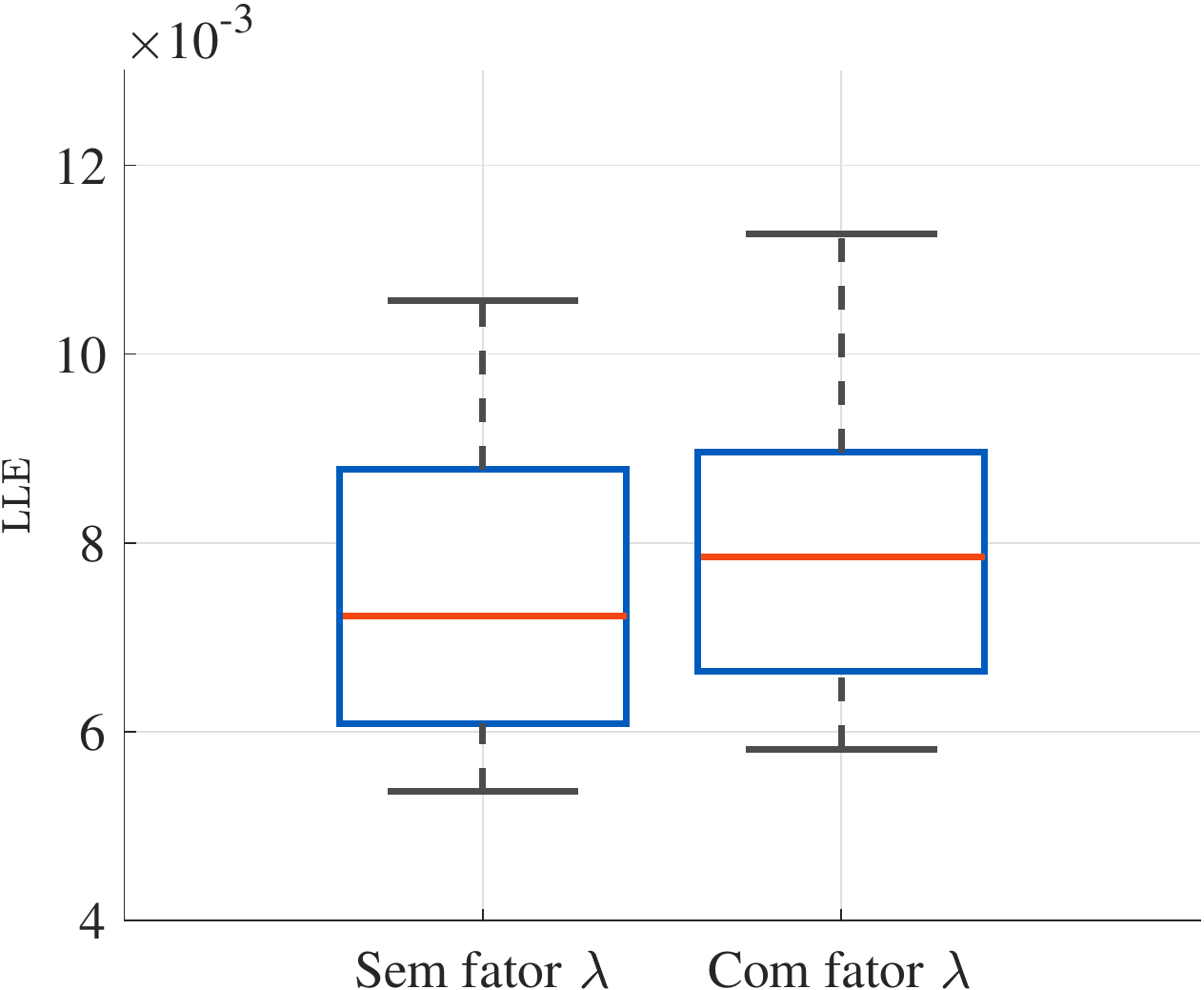} & 
 			\includegraphics[width=0.28\linewidth]{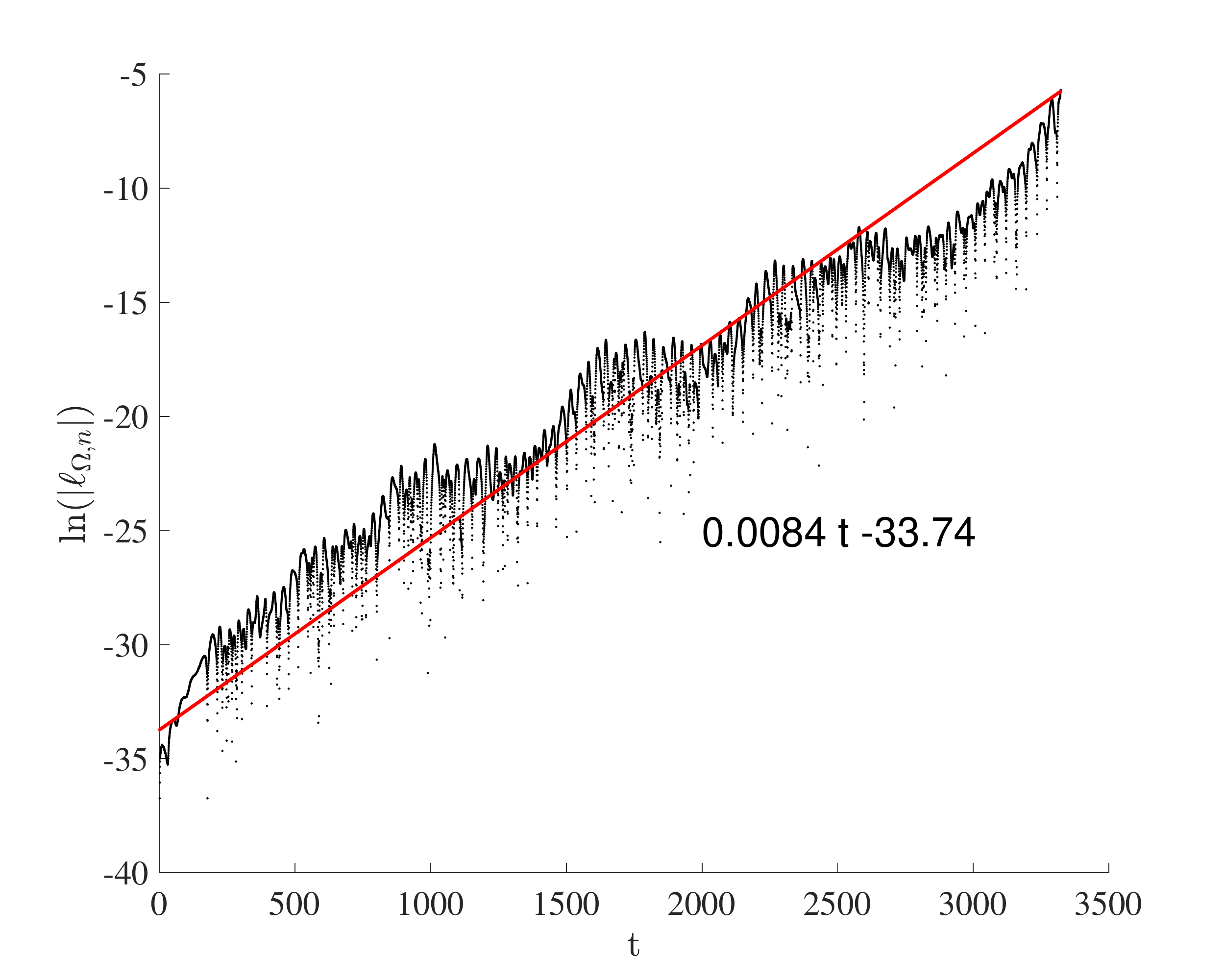}\\
		\end{tabular}        
	\end{center}
\vspace{-0.3cm}
	\caption{Representação gráfica dos valores referentes à média e desvio padrão para os valores de LLE para os métodos sem e com o fator $\lambda$, onde (a) Mapa Logístico, (b) Henon, (c) Mapa Seno, (d) Tenda e (e) Mackey- Glass. Os gráficos (f)-(j) mostram o LBE para os sistemas caóticos adotados. A linha vermelha é o ajuste de mínimos quadrados. Na equação da linha, o primeiro termo é a estimativa do LLE.} 
	\label{fig:media}
\end{figure*}

É possível notar que os valores médios encontrados estão próximos ao esperado na literatura. O desvio padrão fornece a informação da homogeneidade dos dados, é uma medida que expressa o grau de dispersão, indicando, portanto o quanto o conjunto de dados é uniforme. Os valores calculados do desvio padrão são menores quando o fator $\lambda$ é empregado, isso indica que para as diferentes condições iniciais o LLE é mais uniforme quando os dados têm pesos diferentes, como mostrados na Figura \ref{fig:media}. A Figura  \ref{fig:media}(f)-(j) mostra o logaritmo natural do LBE e a linha ajustada para obter a inclinação que é o necessário para quantificar a caoticidade dos sistemas. Nesses gráficos foram usadas as mesmas condições iniciais utilizadas por \cite{PEIXOTO201836}.

\section{Conclusões}

Investigou-se nesse trabalho a ponderação diferenciada dos dados  no cálculo do LLE utilizando modos de arredondamento e estimador recursivos de mínimos quadrados. Cinco sistemas foram testados. Em todos eles o método proposto obteve resultados compatíveis com os da literatura. Além disso, quando comparado com os resultados sem o fator de esquecimento variável, a técnica proposta mostrou-se mais precisa em três sistemas, com pequena variação nos outros dois. Por outro lado, observou-se um aumento da variância em alguns casos evidenciada na Fig. \ref{fig:media}. Importante destacar que neste trabalho também houve a preocupação  de se avaliar o efeito da condição inicial no cálculo do LLE. Isso permite uma rigorosa comparação dos métodos do ponto de vista estatístico. 

Em trabalhos futuros pretende-se avaliar mecanismos para reduzir o número de dados necessários para o cálculo do LLE. Outra questão que merece reflexão é a possibilidade do ajuste não-linear dos dados apresentados na Fig. \ref{fig:media}(f)-(j). Isso poderia gerar valores distintos de LLE para diferentes regiões do espaço de estados.


\end{document}